\documentclass[11pt]{amsart}
\usepackage{amsmath, amssymb, amsfonts, amsthm, bm, url}

\newcommand{\fhi}{\varphi}

\newcommand{\norm}[1]{\lVert#1\rVert}
\newcommand{\ipr}[2]{\left\langle #1, #2 \right\rangle}

\newcommand{\numbersystem}[1]{\mathbb{#1}}
\newcommand{\Quat}{\numbersystem{H}}
\newcommand{\R}{\numbersystem{R}}

\newcommand{\C}{\numbersystem{C}}

\newcommand{\abs}[1]{\left\lvert#1\right\rvert}

\newcommand{\myangle}{\sphericalangle}

\newcommand{\card}[1]{\lvert#1\rvert}
\newcommand{\conj}[1]{\overline{#1}}

\theoremstyle{plain}
\newtheorem{theorem}{Theorem}

\theoremstyle{definition}

\begin{document}

\bibliographystyle{amsplain}

\title[Elementary incidence theorems]{Elementary incidence theorems for complex numbers and quaternions}
\author{J\'ozsef Solymosi}
\thanks{J.~Solymosi was partially supported by
NSERC and OTKA grants and by a Sloan Fellowship.}
\address{Department of Mathematics,
        The University of British Columbia,
        1984 Mathematics Road,
        Vancouver, BC, Canada V6T 1Z2}
\email{\texttt{solymosi@math.ubc.ca}}
\author{Konrad J. Swanepoel}
\thanks{K.~J.~Swanepoel was supported by the South African National Research Foundation under Grant number 2053752.
He also gratefully acknowledges the invitation to and support at the Department of Mathematics, University of British Columbia.}
\address{Department of Mathematical Sciences,
        University of South Africa, PO Box 392,
        Pretoria 0003, South Africa}
\email{\texttt{swanekj@unisa.ac.za}}

\begin{abstract}
We present some elementary ideas to prove the following Sylvester-Gallai type theorems involving incidences between points and lines in the planes over the complex numbers and quaternions.
\begin{enumerate}
\item Let $A$ and $B$ be finite sets of at least two complex numbers each.
Then there exists a line $\ell$ in the complex affine plane such that  $\card{(A\times B)\cap\ell}=2$.
\item Let $S$ be a finite noncollinear set of points in the complex affine plane.
Then there exists a line $\ell$ such that $2\leq\card{S\cap\ell}\leq 5$.
\item Let $A$ and $B$ be finite sets of at least two quaternions each.
Then there exists a line $\ell$ in the quaternionic affine plane such that $2\leq\card{(A\times B)\cap\ell}\leq 5$.
\item Let $S$ be a finite noncollinear set of points in the quaternionic affine plane.
Then there exists a line $\ell$ such that $2\leq\card{S\cap\ell}\leq 24$.
\end{enumerate}
\end{abstract}

\maketitle

\section{Introduction}
The Sylvester-Gallai theorem \cite{Sylvester, Steinberg} states that for any finite noncollinear set of points in the plane, there is a line intersecting the set in exactly two points.
See \cite[\S7.2]{BMP} for a recent list of references regarding this theorem and its relatives.
It is interesting to ask what happens if the underlying field is not the real numbers.
In this paper we consider the complex numbers $\C$ and the quaternions $\Quat$.
For more general considerations see \cite{PS}.
It is a consequence of an inequality of Hirzebruch \cite{Hirzebruch} that any finite noncollinear subset of $\C^2$ always intersects some line in $2$ or $3$ points.
No elementary proof is known of either the inequality or this consequence (see \cite[\S7.3, Problem 7]{BMP}).
However, by adapting Kelly's proof \cite{Coxeter} of the Sylvester-Gallai theorem that uses the smallest perpendicular distance between points and non-incident lines, it is possible to prove the following somewhat weaker statement.

\begin{theorem}\label{complex}
Let $S$ be a finite noncollinear subset of $\C^2$.
Then there exists a line $\ell$ such that 
$2\leq\card{S\cap\ell}\leq 5$.
\end{theorem}

Using the above consequence of the inequality of Hirzebruch, Kelly \cite{Kelly} showed that for any finite nonplanar subset of $\C^3$ there exists a line intersecting the set in exactly two points, thus answering a question of Serre \cite{Serre}.
Another proof can be found in \cite{EPS} that avoids the use of Hirzebruch's inequality.
There it is also shown that in a finite subset of $\Quat^4$ that does not lie on a hyperplane there is a line intersecting the set in exactly two points.
However, it is not known what can be said about $3$-dimensional or $2$-dimensional finite subsets in quaternionic space.
The following theorem is a first step.
It also uses Kelly's shortest distance idea.
\begin{theorem}\label{quaternion}
Let $S$ be a finite noncollinear subset of $\Quat^2$.
Then there exists a line $\ell$ such that $2\leq\card{S\cap\ell}\leq 24$.
\end{theorem}
Note that although the number $24$ is most likely too high, no example is known that shows it has to be larger than $3$.

As easier special cases one may consider grids, i.e., Cartesian products.
Here there are much stronger conclusions.
\begin{theorem}\label{complex-grid}
Let $A, B\subset\C$ with $2\leq\card{A}, \card{B}<\infty$.
Then there exists a line $\ell$ in $\C^2$ such that
$\card{(A\times B)\cap\ell}=2$.
\end{theorem}

\begin{theorem}\label{quaternion-grid}
Let $A, B\subset\Quat$ with $2\leq\card{A}, \card{B}<\infty$.
Then there exists a line $\ell$ in $\Quat^2$ such that
$2\leq\card{(A\times B)\cap\ell}\leq 5$.
\end{theorem}

\section{Proofs}
We consider $\C^2$ as a complex inner product space in the usual way, i.e., we define the inner product of vectors $(a_1,b_1)$ and $(a_2,b_2)$ in $\C^2$ as
\[ \ipr{(a_1,b_1)}{(a_2,b_2)} = a_1\conj{a_2}+b_1\conj{b_2},\]
the norm of $(a,b)\in\C^2$ as
\[ \norm{(a,b)} = \sqrt{\ipr{(a,b)}{(a,b)}}, \]
and the distance between points $(a_1,b_1)$ and $(a_2,b_2)$ in $\C^2$ as $\norm{(a_1-a_2,b_1-b_2)}$.
Working carefully with inner products one can then determine the shortest distance between a point and a line in $\C^2$, as will be done in the next proof.

\begin{proof}[Proof of Theorem~\ref{complex}]
As in Kelly's proof of the real Sylvester-Gallai theorem, we choose three noncollinear points $p,q,r\in S$ such that the distance between $p$ and the line $qr$ is a minimum.
We show that $qr$ contains at most $5$ points.

After a unitary transformation and a dilatation (thus preserving the minimum distance) we may assume that $p=(0,1)$ and that the line $qr$ is the $x$-axis, i.e., $qr\cap S=\{(z_1,0),\dots,(z_k,0)\}$.
Thus the distance between $p$ and $qr$ equals $1$.
Choose any two points $(z_i,0),(z_j,0)\in qr\cap S$.
Then the distance between $(z_j,0)$ and the line $\ell$ through $(0,1)$ and $(z_i,0)$ must be at least $1$.
Considering an arbitrary point
\[(1-\lambda)(z_i,0)+\lambda(0,1)\] on $\ell$, a simple calculation gives the square of its distance to $(z_j,0)$ to be
\begin{align*}
& \norm{(1-\lambda)(z_i,0)+\lambda(0,1)-(z_j,0)}^2\\
=\;& \conj{(z_i-z_j)}^2-(z_i-z_j)\conj{\lambda z_i}-\lambda z_i\conj{(z_i-z_j)}+\abs{\lambda}^2\abs{z_i}^2+\abs{\lambda}^2\\
=\;& \frac{\abs{z_i-z_j}^2}{1+\abs{z_i}^2} + \frac{1+\abs{z_i}^2}{\abs{z_i}^2}\abs{\frac{\abs{z_i}^2}{1+\abs{z_i}^2}(z_i-z_j)-\lambda z_i}^2.
\end{align*}
The last expression is minimised when
\[ \lambda=\frac{(z_i-z_j)\conj{z_i}}{1+\abs{z_i}^2},\]
which gives that the square of the distance between $(z_j,0)$ and $\ell$ is
\[ \frac{\abs{z_i-z_j}^2}{1+\abs{z_i}^2}.\]
Since this quantity is at least $1$ for all distinct $i$ and $j$, we obtain in particular that $\abs{z_i-z_j} > \abs{z_i}$ and $\abs{z_i-z_j} > \abs{z_j}$.
Therefore, for each triangle $\triangle0z_iz_j$ in the Argand plane, the side $z_iz_j$ is the unique longest side.
This implies that the angle $\myangle z_i0z_j > 60^\circ$, giving $\card{qr\cap S}<6$.
\end{proof}

\begin{proof}[Proof of Theorem~\ref{quaternion}]
We consider $\Quat^2$ to be a left vector space, which can be turned into an inner product space in exactly the same way as $\C^2$, where we recall that the conjugate of a quaternion $\alpha=a+bi+cj+dk$ is $\conj{\alpha}=a-bi-cj-dk$.
The distance between two points is defined as before, and we may follow the previous proof almost verbatim, since the commutative law is not used anywhere.
The only change is that the Argand plane has to be replaced by the identification of $\Quat$ with $\R^4$.
We obtain that $\card{qr\cap S}$ is bounded above by the number of rays from the origin in $\R^4$ that are at angles strictly greater than $60^\circ$.
This is easily seen to be bounded above by the kissing number of a $4$-dimensional ball, which was shown to be $24$ by Musin \cite{Musin}.
\end{proof}

We in fact obtained in the above proof that the \emph{strict kissing number} is an upper bound, i.e., the largest number of unit balls in $\R^4$ touching a single unit ball, and none of them touching each other.
Most likely the strict kissing number of $\R^4$ is strictly smaller than $24$, but this is an open problem.

\begin{proof}[Proof of Theorem~\ref{complex-grid}]
Suppose that the conclusion is false.
Then for any points $(a_1,b_1), (a_2,b_2)\in A\times B$ with $a_1\neq a_2$, $b_1\neq b_2$, there exists a third point $(a_3,b_3)\in A\times B$ collinear with $(a_1,b_1)$ and $(a_2,b_2)$.

Colour the elements of $A$ red and the elements of $B$ blue (and to be thought of as points in the Argand plane).
Our assumption easily implies that for any distinct red points $a_1$ and $a_2$ and any distinct blue points $b_1$ and $b_2$ there exists a third red point $a_3$ and blue point $b_3$ such that some orientation preserving similarity $z\mapsto \alpha z+\beta$, $\alpha,\beta\in \C$, maps $a_i$ to $b_i$, $i=1,2,3$.

Choose $a_1$ and $a_2$ to be the pair of red points closest together, and $b_1$ and $b_2$ the pair of blue points furthest apart.
We now analyse the positions of $a_3$ and $b_3$.
Since $a_1$ and $a_2$ are the closest pair of red points, $a_3$ is not in the interior of the circles of radius $a_1a_2$ and with $a_1$ and $a_2$ as centres.
Similarly, $b_3$ is not exterior to to the circles of radius $b_1b_2$ and with centres $b_1$ and $b_2$.
Since there is an orientation-preserving similarity taking $\triangle a_1a_2a_3$ to $\triangle b_1b_2b_3$, there are exactly two possibilities for $a_3$, namely the intersection points of the two circles of radius $a_1a_2$, both making $\triangle a_1a_2a_3$ equilateral.
Similarly, $b_3$ is an intersection point of the two circles of radius $b_1b_2$, and $\triangle b_1b_2b_3$ is equilateral.

Now interchange $b_1$ and $b_2$, i.e., we consider the line through $(a_1, b_2)$ and $(a_2,b_1)$.
By assumption there is a third point $(a_3',b_3')$ on the line through $(a_1,b_2)$ and $(a_2,b_1)$.
As above, $a_3'$ is one of the two points of intersection of the circles with radius $a_1a_2$, and $b_3'$ one of the points of intersection of the circles with radii $b_1b_2$.

Suppose that $a_3'=a_3$.
There is a unique orientation preserving similarity taking $a_2$ to $b_1$ and $a_1$ to $b_2$.
Since this transformation does not map $a_3$ to $b_3$, it follows that $b_3'\neq b_3$ must be the other point of intersection of the two circles of radius $b_1b_2$.
However, the distance $b_3b_3' > b_1b_2$, contradicting the choice of $b_1$ and $b_2$.
It follows that $a_3'\neq a_3$.
Since the similarity taking $a_2$ to $b_1$ and $a_1$ to $b_2$ does not map $a_3$ to $b_3$, it maps $a_3'$ to $b_3$, i.e., $b_3'=b_3$.

We have found a fourth red point $a_3'$ with $\triangle a_1a_2a_3'$ equilateral.
This argument can be repeated indefinitely with other equilateral triangles of smallest side length in $A$, contradicting the finiteness of $A$.
\end{proof}

Note that we made essential use of the fact that the similarities preserve orientation.
Without this, we would only be able to conclude that there is a line $\ell$ such that $2\leq\card{(A\times B)\cap\ell}\leq 3$.

\begin{proof}[Proof of Theorem~\ref{quaternion-grid}]
We proceed along the same lines as in the previous proof.
Since we consider $\Quat^2$ to be a left vector space, lines are either of the form $\{(x,y)\in\Quat^2: x=0\}$ (vertical lines) or of the form $\{(x,y)\in\Quat^2: y=xm+c\}$ for some $m,c\in\Quat$.

Let $\{a_1,a_2\}$ be a closest pair in $A$ and $\{b_1,b_2\}$ a furthest pair in $B$.
Consider the line $\ell$ through $(a_1,b_1)$ and $(a_2,b_2)$.
Since $a_1\neq a_2$, $\ell$ is of the form $y=xm+c$.
Consider the projections
\[ A' = \{x\in\Quat: (x,y)\in\ell\cap(A\times B)\} \]
and
\[ B' = \{y\in\Quat: (x,y)\in\ell\cap(A\times B)\}. \]
Then the mapping
\[ \fhi:\Quat\to\Quat;\; x\mapsto xm+c \]
maps $A'$ to $B'$.
Identifying $\Quat$ with $\R^4$ we see that $\fhi$ is a similarity, which gives that $B'$ is similar to $A'$.
Because the smallest distance in $A'$ (namely $a_1a_2$) is mapped to the largest distance in $B'$ (namely $b_1b_2$), it follows that $A'$ and $B'$ are both equilateral sets, hence $\card{\ell\cap(A\times B)}=\card{A'}=\card{B'}\leq 5$.
The theorem follows.
\end{proof}

Note that although the similarity $\fhi$ in the above proof is orientation-preserving, we have not been able to exploit this fact.


\begin{thebibliography}{10}

\bibitem{BMP}
P.~Brass, W.~Moser, and J.~Pach, \emph{Research problems in discrete geometry}, Springer, New York,
2005.

\bibitem{Coxeter}
H.~S.~M.~Coxeter, \emph{A problem of collinear points}, Amer.\ Math.\ Monthly \textbf{55} (1948), 26--28.

\bibitem{EPS}
 N.~Elkies, L.~M.~Pretorius, K.~J.~Swanepoel, \emph{Sylvester-Gallai Theorems for Complex Numbers and Quaternions}, Discrete Comput.\ Geom.\ \textbf{35} (2006), 361--373.
 
\bibitem{Hirzebruch}
F.~Hirzebruch, \emph{Arrangements of lines and algebraic surfaces}, Arithmetic
  and geometry, Vol. II, Birkh\"auser Boston, Mass., 1983, pp.~113--140.

\bibitem{Kelly}
L.~M.~Kelly, \emph{A resolution of the {S}ylvester-{G}allai problem of
  {J}.-{P}. {S}erre}, Discrete Comput. Geom. \textbf{1} (1986),
  101--104.

\bibitem{Musin}
O.~R.~Musin, \emph{The kissing number in four dimensions}, preprint, 2005. \url{http://arxiv.org/math.MG/0309430} 

\bibitem{PS}
L.~M.~Pretorius, K.~J.~Swanepoel, \emph{The Sylvester-Gallai theorem, colourings and algebra}, submitted, 2006. \url{http://arxiv.org/math.CO/0606131}

\bibitem{Serre}
J.-P.~Serre, \emph{Problems}, Amer. Math. Monthly \textbf{73} (1966), 89.

\bibitem{Steinberg}
  R.~Steinberg,
  Solution to Problem 4065,
  \emph{Amer.\ Math.\ Monthly} \textbf{51} (1944), 169--171.

\bibitem{Sylvester}
  J.~J.~Sylvester,
  \emph{Educational Times} \textbf{46}, No. 383, 156, March 1, 1893. 

\end{thebibliography}

\end{document}